\documentstyle[12pt,twoside]{article}
\newtheorem{theorem}{\bf Theorem}[section]
\newtheorem{proposition}[theorem]{\bf Proposition}

\newtheorem{corollary}[theorem]{\bf Corollary}

\newtheorem{remark}[theorem]{\bf Remark}

\input{amssym}

\date{}
\begin{document}

\title{{\Large\bf Different types  of  weak amenability for Banach algebras}}

\author{{\normalsize\sc M. J. Mehdipour\footnote{Corresponding author} and A. Rejali}}
\maketitle

{\footnotesize  {\bf Abstract.} In this paper, we introduce and investigate the concepts of cyclically weakly amenable and point amenable. Then, we compare these concepts with the concepts of weakly amenable and cyclically amenable and find the relation between them. For example, we prove that a Banach algebra is weakly amenable if and only if it is both cyclically amenable and cyclically weakly amaneble. In the case where $A$ is commutative, the weak amenability and cyclically weak amenability of $A$ are equivalent. We also show that if $A$ is a Banach algebra with $\Delta(A)\neq\emptyset$, then $A$ is cyclically weakly amenable if and only if $A$ is point amenable and essential. For a unital, commutative Banach algebra $A$, the notions of weakly amenable, cyclically weakly amenable and point amenable coincide. In this case, these are equivalent to the fact that every maximal ideal of $A$ is essential.}
{\footnotetext{ 2020 {\it Mathematics Subject Classification}:
 46H20, 46H25, 47B47.

{\it Keywords}: Weak amenability, cyclically amenable, cyclically weakly amenable, point amenable, continuous point derivation.}}

\section{\normalsize\bf Introduction}

Let $D$ be a linear operator  from Banach algebra  $A$ into  $A^*$. Then $D$ is called a \emph{derivation} if $D(ab)=D(a)\cdot b + a\cdot D(b)$ for all $a, b\in A$. Also, $D$ is said to be an \emph{inner derivation} if there exists $F\in A^*$ such that for every $a\in A$
$$
D(a)=\hbox{ad}_F(a):=F\cdot a-a\cdot F,
$$
where for every $b\in A$
$$
\langle F\cdot a, b\rangle=\langle F, ab\rangle\quad\hbox{and}\quad\langle a\cdot F, b\rangle= \langle F, ba\rangle.
$$
If for every $a, b\in A$  
$$
\langle D(a), b\rangle+\langle D(b), a\rangle=0,
$$
then $D$ is called \emph{cyclic}. 
The space of all continuous (inner, cyclic) derivations from $A$ into $A^*$ is denoted by (${\cal I}_{nn}(A, A^*), {\cal Z}_C( A, A^*)$) ${\cal Z}( A, A^*)$, respectively. A Banach algebra $A$ is called \emph{weakly amenable} if
$$
{\cal Z}(A, A^*)={\cal I}_{nn}(A, A^*).
$$
Gronbaek \cite{gro1} introduced the concept of cyclically amenable and studied some study hereditary properties of it. In fact, he called $A$ \emph{cyclically amenable } if 
$$
{\cal Z}_c(A, A^*)={\cal I}_{nn}(A, A^*).
$$
These definitions motivate us to introduce the notion of cyclically weakly amenable. We call $A$ \emph{cyclically weakly amenable} if 
$$
{\cal Z}(A, A^*)={\cal Z}_c(A, A^*).
$$
Note that the mapping $\delta: A\times A\rightarrow \Bbb{C}$ defined by
$$
\delta(a, b)=\langle D(a), b\rangle+\langle D(b), a\rangle
$$
is bilinear. So there exists the bounded linear functional $\tilde{\delta}$ on  the projective tensor product $A\hat{\otimes}A$ such that $$\tilde{\delta}(a\otimes b)=\langle D(a), b\rangle+\langle D(b), a\rangle.$$ Hence we may define the linear operator $T: {\cal Z}( A, A^*)\rightarrow (A\hat{\otimes}A)^*$ by $T(D)=\tilde{\delta}$. One can prove that 
$$
\hbox{ker}T={\cal Z}_C( A, A^*)
$$
and the spaces $\frac{{\cal Z}( A, A^*)}{{\cal Z}_C( A, A^*)}$ and $\hbox{Im}(T)$ are isomorphic as Banach spaces. This implies that  $A$ is cyclically weakly amenable if $T=0$. Also, $A$ is cyclically amenable if
$$
\hbox{ker}T={\cal I}_{nn}(A, A^*).
$$
For $F_1, F_2\in A^*$, we define the continuous linear operator $\frak{D}_{F_1, F_2}: A\rightarrow A^*$ by $$\langle\frak{D}_{F_1, F_2}(a), b\rangle=F_1(a)F_2(b).$$ 
Let  $\Delta(A)$ be the set of all non-zero multiplicative linear functionals on $A$.
Let us recall that a linear functional $d\in A^*$ is called a \emph{continuous point derivation of} $A$ \emph{at} $\varphi\in\Delta_0(A):=\Delta(A)\cup\{0\}$ if for every $a, b\in A$, we have
$$
d(ab)=d(a)\varphi(b)+\varphi(a) d(b).
$$
We call a Banach algebra $A$ \emph{point amenable} if every continuous point derivation of $A$ at $\varphi\in\Delta(A)$ of $A$ is zero. Also,  $A$ is called 0\emph{-point amenable} if every continuous point derivation of $A$ at  $\varphi\in\Delta_0(A)$ is zero. Clearly, 0-point amenability implies point amenability. It is obvious that if $A$ is essential, then 0-point amenability and point amenability are equivalent. It is well-known that every weakly amenable Banach algebra is point amenable; see for example Theorem 2.8.63 in \cite{d}.

Let  $\frak{A}$ be an algebra and $E$ be an $\frak{A}-$bimodule.  Hochschild \cite{ho1,ho2} studied the space ${\cal H}^n(\frak{A}, E)$, cohomology groups of $\frak{A}$ with coefficients in $E$. He gave conditions under which ${\cal H}^n(\frak{A}, E)$ vanish for all $\frak{A}-$bimodule $E$. Later, some authors investigate ${\cal H}^n(\frak{A}, E)$ for Banach algebras. For example,  Kamowitz \cite{kam} investigated Hochschild's results for commutative Banach algebras; see also \cite{kams}. Johnson \cite{joh} studied this space for the Banach algebra $\frak{A}=L^1(G)$, where $G$ is a locally compact group. He showed that $${\cal H}^1(L^1(G), E^*)$$ vanish for all Banach $L^1(G)-$bimodule $E$ if and only if $G$ is amenable. This result leads to use the word ``amenable" for Banach algebras with the property ${\cal H}^1(\frak{A}, E)$ vanish for all Banach $\frak{A}-$bimodule $E$.  Several authors gave the concepts containing amenability. For instance, Bade, Curtis and Dales \cite{bcd} introduced the concept of weak amenability. They considered commutative Banach algebras and gave a definition for weakly amenable Banach algebras. Johnson \cite{j1} generalized this concept for arbitrary Banach algebras. Many authors studied the weak amenability of Banach algebras and investigated the obtain results in the theory amenable Banach algebras for weakly amenable Banach algebras \cite{gro3, gro}; see \cite{gro2} for weak amenability noncommutative Banach algebra, see \cite{gro1} and \cite{s2, s3, s4} for weak amenability of algebras of nuclear operators and Lipschitz algebra, respectively; \cite{fm} for triangular Banach algebra; see also \cite{f, fr, fss} and \cite{b, gro3, p2} for weak amenability of the Fourier algebra and Banach algebra related to groups and semigroups, respectively.

This paper is organized as follow. In Section 2, we study cyclically weakly amenable and point amenable Banach algebras. We show that every cyclically weakly amenable Banach algebra is essential. We also characterize continuous point derivations and  prove that Banach algebra $A$ is point amenable if and only if ${\cal Z}(A, \Bbb{C}_\varphi)=\{0\}$ for all $\varphi\in\Delta(A)$. In Section 3, we find the relation between cyclically weakly amenable and point amenable Banach algebras. For a Banach algebra $A$ with $\Delta(A)\neq\emptyset$, we show that $A$ is cyclically weakly amenable if and only if $A$ is point amenable and essential. In Section 4, we compare the concepts weakly amenable, cyclically amenable, cyclically weakly amenable and point amenable. We establish that $A$ is weakly amenable if and only if $A$ is cyclically amenable and cyclically weakly amenable. In the case where $A$ is commutative, weak amenability and cyclic weak amenability of $A$ coincide. We also verify that if $A$ is unital commutative, then $A$ is weakly amenable if and only if $A$ is cyclically weakly amenable; or equivalently, every maximal ideal of $A$ is essential. These facts are necessary and sufficient conditions that  $A$ be point amenable. Section 5 discusses type of amenability of certain Banach algebras.

\section{\normalsize\bf Cyclically weakly amenable and point amenable}

 A continuous derivation $D: A\rightarrow A^*$ \emph{vanishes on the diameter of} $A$ if
$\langle D(a), a\rangle=0$ for all $a\in A$.

\begin{proposition}\label{cd} Let $A$ be a Banach algebra and $D: A\rightarrow A^*$ be a continuous derivation. Then the following assertion are equivalent.

\emph{(a)} $D$ is cyclic.

\emph{(b)} $D$ vanishes on the diameter of $A$.\\
If $A$ is unital, then the statements \emph{(a)} and \emph{(b)} are equivalent to

\emph{(c)} $\langle D(a), 1\rangle=0$ for all $a\in A$.
\end{proposition}
{\it Proof.} It is clear that (a) implies (b). Let (b) hold. Hence for every $a, b\in A$, we have
\begin{eqnarray*}
\langle D(a), b\rangle+\langle D(b), a\rangle&=&
\langle D(a), a\rangle+\langle D(a), b \rangle\\
&+&\langle D(b), a\rangle+\langle D(b), b\rangle\\
&=&\langle D(a+b), a+b\rangle\\
&=&0.
\end{eqnarray*}
This shows that $D$ is cyclic. Thus (b) implies (a).

Now, let $A$ be unital. If $D$ is cyclic, then for every $a\in A$
$$
\langle D(a), 1\rangle=\langle D(a), 1\rangle+\langle D(1), a\rangle=0.
$$
So (a) implies (c). Finally, let (c) hold. Then
\begin{eqnarray*}
2\langle D(a), a\rangle&=&\langle 2D(a)a, 1\rangle\\
&=&\langle D(a^2), 1\rangle\\
&=&0
\end{eqnarray*}
for all $a\in A$. Therefore, $D$ is cyclic. That is, (c) implies (b).$\hfill\square$\\

An immediate consequence of  Proposition \ref{cd} is the following.

\begin{corollary} \label{asr} Let $A$ be a Banach algebra. Then the following assertions are equivalent.

\emph{(a)} $A$ is  cyclically weakly amenable

\emph{(b)} Every continuous derivation from $A$ into $A^*$ vanishes on the diameter of $A$.\\
If $A$ is unital, then the statements \emph{(a)} and \emph{(b)} are equivalent to

\emph{(c)} $\langle D(a), 1\rangle=0$ for all continuous derivations $D: A\rightarrow A^*$ and $a\in A$.
\end{corollary}

Let us remark that if $A$ is weakly amenable, then $A$ is essential, i.e. $A^2$ is dense in $A$; see for example Theorem 2.8.63 in \cite{d}. The converse holds for commutativity Banach algebras \cite{gro, gro4}. Now, we prove this result for  cyclically weakly amenable Banach algebras.

\begin{proposition}\label{es} Let $A$ be a Banach algebra. Then the following statements hold.

\emph{(i)}  If $A$ is  cyclically weakly amenable, then $A$ is essential.

\emph{(ii)}  If $A^{**}$ is  cyclically weakly amenable, then $A$ is essential.
\end{proposition}
{\it Proof.} Suppose that $A$ is not essential. Choose $a_0\in A$ such that $a_0\notin \overline{A^2}$. In view of Hahn-Banach theorem, there exists $F\in A^*$ such that $F(a_0)=1$ and $F(A^2)=\{0\}$. Consider the linear operator $\frak{D}_{F, F}: A\rightarrow A^*$. Since $F|_{A^2}=0$, for every $a, b, c\in A$, we have
$$
\frak{D}_{F, F}(ab)=F(ab)F=0
$$
and
\begin{eqnarray*}
\langle a\frak{D}_{F, F}(b), c\rangle+\langle \frak{D}_{F, F}(a)b, c\rangle&=&\langle \frak{D}_{F, F}(b), ca\rangle+\langle \frak{D}_{F, F}(a), bc\rangle\\
&=& F(b)F(ca)+F(a)F(bc)\\
&=&0.
\end{eqnarray*}
Hence $\frak{D}_{F, F}$ is a continuous derivation. Thus $\frak{D}_{F, F}$ is cyclic. It follows that
\begin{eqnarray*}
\langle \frak{D}_{F, F}(a_0), a_0\rangle&+&\langle \frak{D}_{F, F}(a_0) , a_0\rangle\\&=&2\langle \frak{D}_{F, F}(a_0), a_0\rangle\\
&=& 2F(a_0)F(a_0)=2,
\end{eqnarray*}
a contradiction. This shows that (i) holds. The argument used in the proof of Proposition 2.1 in \cite{gl}  proves  (ii). So we omit it. $\hfill\square$\\

Let $\varphi$ be a character of Banach algebra $A$ with kernel of $M_\varphi$. It is easy to see that $\Bbb{C}$ with the following actions is a Banach $A-$bimodule which is denoted by $\Bbb{C}_\varphi$.
$$
a\cdot\alpha=\alpha\cdot a=\varphi(a)\alpha
$$
for all $a\in A$ and $\alpha\in\Bbb{C}$.

\begin{proposition}\label{pit} Let $A$ be a Banach algebra, $\varphi\in\Delta(A)$ and $d\in A^*$. Then the following assertions are equivalent.

\emph{(a)}  $d$ is a continuous point derivation at  $\varphi$ of $A$.

\emph{(b)} $d$ is a continuous derivation from $A$ into $\Bbb{C}_\varphi$.

\emph{(c)}  $\frak{D}_{d, \varphi}$ is a continuous derivation.\\
 In this case, $d|_{M_{\varphi}^2}=0$.

\end{proposition}
{\it Proof.}  It is easy to see that a continuous linear functional $d$ on $A$ is a continuous point derivation at $\varphi$ if and only if $d\in{\cal Z}(A, \Bbb{C}_\varphi)$. That is, (a) and (b) are equivalent. For every $a, b, x\in A$, we have
\begin{eqnarray}\label{neko}
\langle \frak{D}_{d, \varphi}(a)b+a\frak{D}_{d, \varphi}(b), x\rangle&=& d(a)\varphi(bx)+ d(b)\varphi(xa)\\
&=&(d(a)\varphi(b)+d(b)\varphi(a))\varphi(x).\nonumber
\end{eqnarray}
Let $\frak{D}_{d, \varphi}$ be a continuous derivation. Choose $x\in A$ with $\varphi(x)=1$. Then by (\ref{neko}) we have
\begin{eqnarray*}
d(ab)=\langle \frak{D}_{d, \varphi}(ab), x\rangle=d(a)\varphi(b)+d(b)\varphi(a)
\end{eqnarray*}
for all $a, b\in A$. Hence $d$ is a continuous point derivation at  $\varphi$ of $A$.

Conversely, let $d$ be a continuous point derivation of $A$ at $\varphi$. Then for every $a, b, x\in A$, we have
$$
(d(a)\varphi(b)+d(b)\varphi(a))\varphi(x)=d(ab)\varphi(x)=\langle \frak{D}_{d, \varphi}(ab), x\rangle.
$$
This together with (\ref{neko}) shows that $\frak{D}_{d, \varphi}$ is a continuous derivation.

To complete the proof, let $x, y\in M_\varphi=\hbox{ker}\varphi$. Then $\varphi(x)=\varphi(y)=0$ and so
$$
d(xy)=d(x)\varphi(y)+\varphi(x) d(y)=0.
$$
This implies that $d|_{M_\varphi^2}=0$.$\hfill\square$\\

In the sequel, we investigate Proposition \ref{pit} in some special cases.

\begin{proposition}\label{malmal} Let $A$ be a  Banach algebra, $\varphi\in\Delta(A)$ and $d$ be a non-zero linear functional on $A$. Consider the following assertions.

\emph{(a)}  $d$ is a continuous point derivation at  $\varphi$ of $A$.


\emph{(b)} $\frak{D}_{d, \varphi}$ is a non-inner continuous derivation.

\emph{(c)} $\frak{D}_{d, \varphi}$ is a non-cyclic continuous derivation.\\
\emph{(i)}  If $A$ is commutative, then the statements \emph{(a)} and \emph{(b)} are equivalent.\\
\emph{(ii)} If $A$ is essential, then the statements \emph{(a)} through \emph{(c)} are equivalent.\\
\emph{(iii)} If $A$ is unital, then the statements \emph{(a)} through \emph{(c)} are equivalent to

\emph{(d)}  $d(1_A)=0$ and  $d|_{M_{\varphi}^2}=0$.
\end{proposition}
{\it Proof.} The implications (b)$\Rightarrow$ (a) and (c)$\Rightarrow$ (a) follow from Proposition \ref{pit}. Now, let $d$ be a continuous point derivation at  $\varphi$ of $A$ and $\frak{D}_{d, \varphi}$  be an inner derivation. Then $\frak{D}_{d, \varphi}=\hbox{ad}_F$ for some $F\in A^*$. So for every $a, b\in A$, we have
$$
d(a)\varphi(b)=\langle \frak{D}_{d, \varphi}(a), b\rangle=\langle F\cdot a-a\cdot F, b\rangle=\langle F, ab-ba\rangle.
$$
This shows that if $A$ is unital or commutative, then $d=0$, a contradiction. That is, (a)$\Rightarrow$(b). 

Assume that $A$ is unital. If $d$ is a continuous point derivation at $\varphi\in\Delta(A)$, then
$$
d(1_A)=d(1_A1_A)=2d(1_A)\varphi(1_A)=2d(1_A).
$$
So $d(1_A)=0$. By Proposition \ref{pit}, $d|_{M_\varphi^2}=0$. So (a)$\Rightarrow$ (d). Note that if $d\in A^*$ such that $d(1_A)=0$ and $d|_{M_\varphi^2}=0$, then $d|_{\Bbb{C}1_A+M_\varphi^2}=0$. Thus $d$ is a continuous point derivation at  $\varphi$ of $A$; see Proposition 1.8.8 in \cite{d}. Thus (d)$\Rightarrow$ (a).

Finally, let $A$ be essential and (a) hold. If  $\frak{D}_{d, \varphi}$ is a cyclic continuous derivation, then for every $a, b\in A$, we have
\begin{eqnarray*}
d(ab)&=& d(a)\varphi(b)+\varphi(a) d(b)\\
&=&\langle\frak{D}_{d, \varphi}(a), b\rangle+\langle\frak{D}_{d, \varphi}(b), a\rangle=0.
\end{eqnarray*}
So $d$ is zero on $A^2$. Since $A$ is essential, $d=0$ on $A$.
This contradiction shows that (a)$\Rightarrow$ (c).$\hfill\square$\\

 As a consequence of Proposition \ref{pit}  we give the following result.

\begin{corollary}\label{mokh} Let $A$ be a Banach algebra with $\Delta(A)\neq\emptyset$. Then the following assertions are equivalent.

\emph{(a)} $A$ is point amenable.

\emph{(b)} ${\cal Z}(A, \Bbb{C}_\varphi)=\{0\}$ for all $\varphi\in\Delta(A)$.

\emph{(c)} $\frak{D}_{d, \varphi}=0$  for all  continuous point derivation d at $\varphi\in\Delta(A)$ of $A$.\\
\end{corollary}

Let $F_i$ be a bounded linear functional on Banach algebra $A_i$, for $i=1,2$. In the following, we define the bounded linear functional $F_1\otimes F_2$ on $A_1\hat{\otimes}A_2$ by $$F_1\otimes F_2(a_1\otimes a_2)=F_1(a_1)F_2(a_2)$$ for all $a_1\in A_1$ and $a_2\in A_2$.

\begin{theorem}\label{mahgom} Let $A_i$ be  a Banach algebra and $d_i$ be a continuous point derivation at $\varphi_i\in\Delta_0(A_i)$, for $i=1, 2$. Then $d_1\otimes\varphi_2+d_2\otimes\varphi_1$ is a continuous point derivation at $\varphi_1\otimes\varphi_2$ of $A_1\hat{\otimes}A_2$.
\end{theorem}
{\it Proof.} Let $d_i$ be a continuous point derivation at $\varphi_i\in\Delta_0(A_i)$, for $i=1, 2$. Set $d:=d_1\otimes\varphi_2+d_2\otimes\varphi_1$ and $\phi:=\varphi_1\otimes\varphi_2$. For every $a_1, x_1\in A_1$ and $a_2, x_2\in A_2$, we have
\begin{eqnarray*}
d((a_1\otimes a_2)(x_1\otimes x_2))&=&
d(( (a_1x_1)\otimes (a_2x_2)))\\
&=&
d_1( a_1x_1)\varphi_2(a_2x_2)+\varphi_1(a_1x_1)d_2(a_2x_2)\\
&=&
\varphi_1(a_1)d_1(x_1)\varphi_2(a_2)\varphi_2(x_2)+
d_1(a_1)\varphi_1(x_1)\varphi_2(a_2)\varphi_2(x_2) \\
&+&
d_2(a_2)\varphi_2(x_2)\varphi_1(a_1)\varphi_1(x_1)+
d_2(x_2)\varphi_2(a_2)\varphi_1(a_1)\varphi_1(x_1).
\end{eqnarray*}
Also,
\begin{eqnarray*}
\phi(a_1\otimes a_2)d( x_1\otimes x_2)&+&
\phi( x_1\otimes x_2)d(a_1\otimes a_2)\\
&=&
\varphi_1(a_1)\varphi_2(a_2)d_1\otimes\varphi_1(x_1\otimes x_2)\\
&+&
\varphi_1(a_1)\varphi_2(a_2)d_2\otimes\varphi_1(x_1\otimes x_2)\\
&+&
\varphi_1(x_1)\varphi_2(x_2)d_1\otimes\varphi_1(a_1\otimes a_2)\\
&+&
\varphi_1(x_1)\varphi_2(x_2)d_2\otimes\varphi_1(a_1\otimes a_2)\\
&=&
\varphi_1(a_1)\varphi_2(a_2)d_1(x_1)\varphi_2(x_2)\\
&+&
\varphi_1(a_1)\varphi_2(a_2)d_2(x_1)\varphi_1(x_2)\\
&+&
\varphi_1(x_1)\varphi_2(x_2)d_1(a_1)\varphi_2(a_2)\\
&+&
\varphi_1(x_1)\varphi_2(x_2)d_2(a_1)\varphi_1(a_2).
\end{eqnarray*}
Therefore, $d$ is a continuous point derivation at $\phi$ of $A_1\hat{\otimes}A_2$.$\hfill\square$

\section{\normalsize\bf Relation between cyclically weakly amenable and point amenable}

Let $A$ be a Banach algebra. We call a continuous linear functional $p\in (A\hat{\otimes}A)^*$  \emph{quasi-additive functional of} $A$ if  for every $a, b, x\in A$,
$$
p(ax\otimes b)=p(a\otimes xb)+p(x\otimes ba).
$$
A quasi-additive functional $p$ of  $A$ is called \emph{inner} if there exists $F\in A^*$ such that $$p(a\otimes b)=F(ab-ba)$$ for all $a, b\in A$. Also, $p$ is called \emph{cyclic} if $p(a\otimes a)=0$ for all $a\in A$.

\begin{theorem}\label{helma} Let $D: A\rightarrow A^*$ be a continuous linear operator, $d\in A^*$ and $\varphi\in\Delta(A)$. Then the following statements hold.

\emph{(i)} $D$ is a derivation if and only if there exists a quasi-additive functional $p$ of $A$ such that $p(a\otimes b)=\langle D(a), b\rangle$ for all $a, b\in A$.

\emph{(ii)} $D$ is an inner derivation if and only if there exists an inner quasi-additive functional $p$ of $A$ such that $p(a\otimes b)=\langle D(a), b\rangle$ for all $a, b\in A$.

\emph{(iii)} $D$ is a cyclic derivation  if and only if there exists a cyclic quasi-additive functional $p$ of $A$ such that  $p(a\otimes b)=\langle D(a), b\rangle$ for all $a, b\in A$.

\emph{(iv)} If $A$ is unital, then $D$ is a cyclic derivation  if and only if there exists a quasi-additive functional $p$ of $A$ such that  $p(a\otimes b)=\langle D(a), b\rangle$ and $p(a\otimes 1)=0$ for all $a, b\in A$.

\emph{(v)} $d$ is a continuous point derivation at $\varphi\in\Delta(A)$ if and only if there exist $a_0\in  A\setminus M_\varphi$ and  a quasi-additive functional $p$ of $A$ such that $d(a)=\frac{p(a\otimes a_0)}{\varphi(a_0)}$ for all $a\in A$.
\end{theorem}
{\it Proof.} (i) It is well-known that the function
$$
\Gamma: {\cal B}(A, A^*)\rightarrow (A\hat{\otimes}A)^*
$$
defined by
$$
\langle \Gamma(T), a\otimes b\rangle=\langle T(a), b\rangle
$$
is an isometric isomorphism as Banach spaces; see for example Proposition 13 VI in \cite{bd}. So if $D$ is a continuous derivation from $A$ into $A^*$, then $D\in{\cal B}(A, A^*)$ and hence 
 $$
\langle \Gamma(D), a\otimes b\rangle=\langle D(a), b\rangle
$$
for all $a, b\in A$. Put $p:=\Gamma(D)$. It is routine to check that $p$ is a quasi-additive functional of $A$. The converse is clear.

(ii) Let $D$ be an inner derivation. Thus there exists $F\in A^*$ such that 
$$
\langle D(a), b\rangle=F(ab-ba)
$$
for all $a, b\in A$. This together with (i) proves (ii).

(iii) This follows from (i) and Corollary \ref{asr}.

(iv) This follows from (i) and Corollary \ref{asr}.

(v) Let $d$ be a continuous point derivation at $\varphi\in\Delta(A)$. Then the linear operator $\frak{ D}_{d, \varphi}$ is a continuous derivation from $A$ into $A^*$. By (i), there exists a quasi-additive functional $p$ of $A$ such that 
$$
p(a\otimes b)=\langle\frak{D}_{d, \varphi}(a), b\rangle=d(a)\varphi(b)
$$
for all $a, b\in A$. Choose $a_0\in  A\setminus M_\varphi$. Then $\varphi(a_0)$ is non-zero and
$$d(a)=\frac{p(a\otimes a_0)}{\varphi(a_0)}$$ for all $a\in A$. The converse is clear.$\hfill\square$\\

The next result is an immediate consequence of Theorem \ref{helma}.

\begin{corollary}\label{asrs} Let $A$ be a Banach algebra. Then the following statements hold.

\emph{(i)} $A$ is weakly amenable if and only if every quasi-additive functional of $A$ is inner.

\emph{(ii)} $A$ is cyclically weakly amenable if and only if every quasi-additive functional of $A$ is cyclic.

 \emph{(iii)} $A$ is cyclically amenable if and only if every cyclic quasi-additive functional of $A$ is inner.

\emph{(iv)} In the case where $\Delta(A)\neq\emptyset$, $A$ is point amenable if and only if for every $\varphi\in\Delta(A)$, $b\in  A\setminus M_\varphi$ and quasi-additive functional $p$ of $A$, the function $a\mapsto p(a\otimes b)$ is zero on $A$.
\end{corollary} 

\begin{theorem}\label{mardom1} Let $A$ be a Banach algebra with $\Delta(A)\neq\emptyset$ . Then the following assertions are equivalent.

\emph{(a)} $A$ is cyclically weakly amenable.

\emph{(b)} $A^\sharp$ is cyclically weakly amenable.

\emph{(c)} $A^\sharp$ is point amenable

\emph{(d)} $A^\sharp$ is 0-point amenable.

\emph{(e)} $A$ is 0-point amenable.

\emph{(f)} $A$ is point amenable and essential.
\end{theorem} 
{\it Proof.} Assume that $A$ is a cyclically weakly amenable Banach algebra. By Proposition \ref{es}, $A$ is essential. Now, let $\varphi\in\Delta(A)$. If $d$ is a continuous point derivation at $\varphi$ of $A$, then $\frak{D}_{d, \varphi}$  is a continuous derivation from $A$ into $A^*$. So  $\frak{D}_{d, \varphi}$ is cyclic. From Proposition \ref{malmal}  we see that $\frak{D}_{d, \varphi}=0$ and hence $d=0$. Thus (a)$\Rightarrow$(f). It is obvious that (f)$\Rightarrow$(e) and (d)$\Rightarrow$(c). 

Let $d$ be a continuous point derivation at $\varphi\in\Delta_0(A^\sharp)$ of $A^\sharp$. Define $\bar{d}, \bar{\varphi}: A\rightarrow \Bbb{C}$ by 
$$
\bar{d}(a)=d(a, 0)\quad\hbox{and}\quad\bar{\varphi}(a)=\varphi(a, 0)
$$
for all $a\in A$. It is clear that $\bar{d}$ is a continuous point derivation at $\bar{\varphi}\in\Delta_0(A)$ of $A$. So if (e) holds, then $\bar{d}=0$. Hence $d(a, 0)=0$ for all $a\in A$. But, $d(0,1)=0$. Thus $d=0$ and therefore, $A^\sharp$ is 0-point amenable. That is, (e)$\Rightarrow$(d).

Let $\psi\in\Delta(A^\sharp)$. Then $\psi(1_{A^\sharp})=1$. So $1_{A^\sharp}\in  A^\sharp\setminus M_\psi$. If $A^\sharp$ is point amenable and $p$ is a quasi-additive functional of $A^\sharp$, then by Corollary \ref{asrs} (iv),
$p(a\otimes 1)=0$ for all $a\in A^\sharp$. On the other hand, 
$$
2p(a\otimes a)=p(a^2\otimes 1_{A^\sharp})=0
$$
for all $a\in A^\sharp$. Hence $p$ is cyclic and thus $A^\sharp$ is cyclically weakly amenable. That is, (c)$\Rightarrow$(b). Finally, assume that $A^\sharp$ is cyclically weakly amenable and $D: A\rightarrow A^*$ be a continuous derivation. Then $\bar{D}: A^\sharp\rightarrow (A^\sharp)^*$ define by 
$$
\langle\bar{D}((a, \alpha), (b, \beta)\rangle=\langle D(a), b\rangle\quad(a, b\in A, \alpha, \beta\in\Bbb{C})
$$
is a continuous derivation. Thus $\bar{D}=0$. It follows that $D=0$. So (b)$\Rightarrow$(a).$\hfill\square$\\

Now, we give a consequences of Theorem \ref{mardom1}.

\begin{corollary}\label{bil} Let $A$ be an essential Banach algebra with $\Delta(A)\neq\emptyset$. Then the following assertions are equivalent.

\emph{(a)} $A$ is cyclically weakly amenable.

\emph{(b)}  $A$ is 0-point amenable.

\emph{(c)} $A$ is point amenable
\end{corollary}

For $\varphi\in\Delta(A)$, we denote by $(M_\varphi\setminus M_\varphi^2)^\times$ the set of all linear functionals on $M_\varphi\setminus M_\varphi^2$.

\begin{theorem}\label{sal}  Let $A$ be a unital Banach algebra with $\Delta(A)\neq\emptyset$. Then the following assertions are equivalent.

\emph{(a)} $A$ is cyclically weakly amenable.

\emph{(b)} $A$ is point amenable. 

\emph{(c)} $(M_\varphi\setminus M_\varphi^2)^\times=\{0\}$ for all $\varphi\in\Delta(A)$.
\end{theorem}
{\it Proof.} By Theorem \ref{mardom1} the statements (a) and (b) are equivalent. From Proposition 1.8.8 in \cite{d} we see that ${\cal Z}(A, \Bbb{C}_\varphi)$ is linearly isomorphic to $(M_\varphi\setminus M_\varphi^2)^\times$. This fact together with Corollary \ref{mokh} shows that $A$ is point amenable if and only if $(M_\varphi\setminus M_\varphi^2)^\times=\{0\}$ for all $\varphi\in\Delta(A)$. $\hfill\square$

\section{\normalsize\bf Comparison of type of amenability}

In the following we study the relation between the concepts  weak amenability, weak cyclic amenability and point amenability.

\begin{theorem}\label{shah} Let $A$ be a Banach algebra. Then the following assertions are equivalent.

\emph{(a)} $A$ is weakly amenable.

\emph{(b)}  $A$ is both cyclically amenable and  cyclically weakly amenable.\\
In this case, $A$ is point amenable and essential.
\end{theorem}
{\it Proof.} Let $A$ be weakly amenable. Clearly,  $A$ is cyclically amenable. Let $p$ be a quasi-additive functional of $A$. By Corollary \ref{asrs} (i), there exists $F\in A^*$ such that 
$$
p(a\otimes b)=F(ab-ba)
$$
for all $a, b\in A$. So $p(a\otimes a)=0$. Thus $A$ is cyclically weakly amenable. Therefore, (a) implies (b). An easy application of Corollary \ref{asrs} shows that (b) implies (a). Finally, By Theorem \ref{mardom1} $A$ is point amenable and essential.$\hfill\square$\\

In the following, we show that for commutative Banach algebras the weak amenability and weak cyclic amenability coincide.

\begin{theorem}\label{40}  Let $A$ be a commutative Banach algebra. Then $A$ is weakly amenable if and only if $A$ is cyclically weakly amenable.
\end{theorem}
{\it Proof.}  Let $A$ be  cyclically weakly amenable. Suppose that $D: A\rightarrow A^*$ is non-zero derivation. Since $A$ is essential, there exist $a_0, b_0\in A$ such that $D(a_0b_0)$ is non-zero. Fix $\Phi\in A^{**}$ such that
$$
\langle\Phi, D(a_0b_0)\rangle=1.
$$
Let us remark that $A^*$ is a Banach $A-$bimodule with the following action.
$$
\langle F\cdot a, b\rangle=\langle F, ab\rangle\quad\quad(F\in A^*, a, b\in A).
$$
Then for every $a, b, c, x\in A$, we have
\begin{eqnarray*}
\langle D(ab)\cdot c, x\rangle&=&\langle D(ab), cx\rangle\\
&=&\langle D(a)\cdot b+a\cdot D(b), cx\rangle\\
&=&\langle D(a), b(cx)\rangle+\langle D(b), (cx)a\rangle\\
&=&\langle D(a), b(cx)\rangle+\langle D(b), (ca)x\rangle\\
&=&\langle D(a)\cdot(bc), x\rangle+\langle D(b)\cdot(ca), x\rangle\\
&=&\langle D(a)\cdot(bc)+ D(b)\cdot (ca), x\rangle.
\end{eqnarray*}
This shows that
\begin{eqnarray*}\label{1}
D(ab)\cdot c=D(a)\cdot(bc)+ D(b)\cdot (ca).
\end{eqnarray*}
We define the bounded linear operator $T: A^*\rightarrow A^*$ by
$T(F)=\Phi\vartriangle F$, where $$\Phi\vartriangle F(a)=\langle \Phi, F\cdot a\rangle$$
for all $a\in A$ and $F\in A^*$.
Set $U:=T\circ D$. Then
\begin{eqnarray*}
U(a)=T(D(a))=\Phi\vartriangle D(a).
\end{eqnarray*}
Now for every $a, b, c\in A$ we have
\begin{eqnarray*}
\langle U(a)\cdot b+a\cdot U(b), c \rangle&=&\langle U(a), bc\rangle+\langle U(b), ca \rangle\\
&=& \langle \Phi\vartriangle D(a), bc\rangle+ \langle \Phi\vartriangle D(b), ca\rangle\\
&=&\langle \Phi, D(a)\cdot(bc)+D(b)\cdot (ca)\rangle\\
&=&\langle \Phi, D(ab)\cdot c\rangle\\
&=&\langle \Phi\vartriangle D(ab), c\rangle\\
&=&\langle U(ab), c\rangle.
\end{eqnarray*}
It follows that
$$
U(ab)=U(a)\cdot b+a\cdot U(b).
$$
Therefore, $U$ is a continuous derivation from $A$ into $A^*$. Hence $U$ is cyclic. Thus
\begin{eqnarray*}
\langle U(a_0), b_0\rangle&+&\langle U(b_0), a_0\rangle\\
&=&\langle \Phi\vartriangle D(a_0), b_0\rangle+\langle \Phi\vartriangle D(b_0), a_0\rangle\\
&=&\langle\Phi, D(a_0)\cdot b_0\rangle+\langle\Phi, D(b_0)\cdot a_0\rangle\\
&=&\langle\Phi, D(a_0)\cdot b_0+D(b_0)\cdot a_0\rangle\\
&=&\langle \Phi, D(a_0b_0)\rangle\\
&=&1.
\end{eqnarray*}
This contradiction shows that $D$ is zero. Hence $A$ is weakly amenable.$\hfill\square$\\

Gronbaek \cite{gro} showed that if $A$ is weakly amenable, commutative Banach algebra and $I$ is a closed ideal in $A$, then $I$ is weakly amenable if and only $I$ is essential. This result and Corollary \ref{40} prove the next result.

\begin{corollary} Let $A$ be a commutative   weakly amenable Banach algebra and $I$ be a closed ideal in $A$. Then the following assertions are equivalent.

\emph{(a)} $I$ is weakly amenable.

\emph{(b)} $I$ is cyclically weakly amenable.

\emph{(c)} $I$ is essential.
\end{corollary}

\begin{remark} {\rm Let us recall that a Banach algebra $A$ is called \emph{contractible} if for every Banach $A-$bimodule $E$, every continuous derivation from $A$ into $E$ is inner. It is proved that contractible Banach algebras are unital; see \cite{hel1}.
 There is a claim that if $A$ is a commutative weakly amenable Banach algebra, then for every Banach $A-$bimodule $E$, the only continuous derivation from $A$ into $E$ is zero; see for example Theorem 2.8.63 in \cite{d}.  It seems that this claim can not true, because otherwise every commutative weakly amenable Banach algebra is contractible. Hence these Banach algebras are unital. Now note that $L^1({\Bbb R})$ is a weakly commutative Banach algebra without identity.}
\end{remark}

\begin{proposition}  Let $A$ be a Banach algebra, $\varphi\in\Delta(A)$ and $d$ be a non-zero linear functional on $A$.  If $A$ is weakly amenable, then $\frak{D}_{d, \varphi}$ is not a continuous derivation.
\end{proposition}
{\it Proof.} Suppose that $\frak{D}_{d, \varphi}$ is a continuous derivation. It follows from Proposition \ref{pit} that $d$ is a continuous point derivation at $\varphi$ of $A$. So Proposition \ref{malmal} shows that $\frak{D}_{d, \varphi}=0$ and hence $d=0$, a contradiction.$\hfill\square$\\

In the next result, we show that weakly amenable, cyclically weakly amenable and point amenable of unital, commutative Banach algebras are equivalent.

\begin{theorem}\label{tasisat} Let $A$ be a unital, commutative Banach algebra. Then the following assertions are equivalent.

\emph{(a)} $A$ is weakly amenable.

\emph{(b)} $A$ is cyclically weakly amenable.

\emph{(c)} $A$ is point amenable.

\emph{(d)} $(M_\varphi\setminus M_\varphi^2)^\times=\{0\}$ for all $\varphi\in\Delta(A)$.

\emph{(e)} For every $\varphi\in\Delta(A)$, the maximal ideal $M_\varphi$ is essential.

\emph{(f)}  Every maximal ideal of $A$ is essential.

\emph{(g)} $H^1(A, E)=\{0\}$ for all finite-dimensional Banach $A-$bimodule $E$.
\end{theorem}
{\it Proof.} In view of Theorem \ref{sal} and Corollary \ref{40}, the statements (a) through (d) are equivalent.  It is well-known that the maximal ideals of $A$ is the kernel of elements of $\Delta(A)$; see Theorem 2, page 79 in \cite{bd}. Thus (e) and  (f) are equivalent.  By  Proposition 1.2 in \cite{bcd}, the statements (c) and (g) are equivalent. Now, let (d) hold. Then $(M_\varphi\setminus M_\varphi^2)^*=\{0\}$. Thus $(M_\varphi\setminus M_\varphi^2)^{**}=\{0\}$. It follows that $M_\varphi\setminus M_\varphi^2=\{0\}$. So $M_\varphi^2=M_\varphi$. This implies that $M_\varphi$ is essential. So (d) implies (e). To complete the proof, let $d$ be a continuous point derivation at $\varphi\in\Delta(A)$ of $A$. If $a\in A$, then $a\in M_\varphi$ for some $\varphi\in\Delta(A)$. Since $M_\varphi$ is a maximal ideal  of $A$, by assumption $M_\varphi$ is essential. Thus $M_\varphi^2=M_\varphi$. This together with Proposition \ref{pit} shows that $d|_{M_\varphi}=0$ and so $d(a)=0$. That is, (e) implies (c).$\hfill\square$

\begin{corollary} Let  $A$ be a unital, commutative Banach algebra. Then the following assertions are equivalent.

\emph{(a)}  For every quasi-additive functional $p$ of $A$, $p(a\otimes b)=0$ for all $a, b\in A$.

\emph{(b)}  For every quasi-additive functional $p$ of $A$ and $\varphi\in\Delta(A)$, $p(a\otimes b)=0$ for all $a\in A$ and $b\in  A\setminus M_\varphi$.

\emph{(c)}  For every quasi-additive functional $p$ of $A$, $p(a\otimes a)=0$ for all $a\in A$.

\emph{(d)} For every quasi-additive functional $p$ of $A$, $p(a\otimes 1)=0$ for all $a\in A$.
\end{corollary}

\begin{theorem}\label{sem} Let $A$ be a semisimple Banach algebra. Then the following assertions are equivalent.

\emph{(a)} $A$ is weakly amenable.

\emph{(b)} $A$ is cyclically weakly amenable.

\emph{(c)} $A$ is 0-point amenable.

\emph{(d)} $A$ is point amenable.

\emph{(e)} $A^\sharp$ is weakly amenable.

\emph{(f)} $A^\sharp$ is cyclically weakly amenable.

\emph{(g)} $A^\sharp$ is 0-point amenable.

\emph{(h)} $A^\sharp$ is point amenable.
\end{theorem} 
{\it Proof.} In view of  Theorems \ref{mardom1} and \ref{shah}, it suffices to show that (d)$\Rightarrow$(a). To this end,
let $A$ be point amenable and $p$ be a quasi-additive functional of $A$. Since $A$ is semisimple, the Gelfand homomorphism is one-to-one. So if $b$ is a non-zero element of $A$, then $b\in  A\setminus M_\varphi$ for some $\varphi\in\Delta(A)$. By Corollary \ref{asrs}, $p(a\otimes b)=0$ for all $a\in A$. This implies that $A$ is weakly amenable.$\hfill\square$

\section{\normalsize\bf Type of amenability of certain Banach algebras}

Haagerup \cite{h1} showed that every $C^*-$algebra is weakly amenable; see \cite{h2} for a simple proof. The following result at once follows from \cite{h1}. Before, we give another proof by using the concept of one-point compactification, let us recall that  $C_0(X)$ is the space of all continuous complex-valued functions $f$ on  locally compact Hausdorff space $X$ that vanish at infity \cite{c}.

\begin{theorem}\label{panda} Let $X$ be a locally compact Hausdorff space. Then $C_0(X)$ is weakly amenable.
\end{theorem}
{\it Proof.} Let $M$ be a maximal ideal of $C(X)$. Then $M=\hbox{ker}(\delta_x)$ for some $x\in X$. Hence $M$ is a $C^*-$algebra and so it has a bounded approximate identity. Thus $M$ is essential. By Theorems \ref{tasisat} and \ref{sem}, $C_0(X)$ is weakly amenable.$\hfill\square$\\

Let us remark that the second dual of $C^*-$algebras are  von Neumann algebras; see for example \cite{d}. As an immediate consequence of Theorem \ref{panda} we can obtain Haagerup's theorem for commutative $C^*-$algebras.

\begin{proposition}\label{ccmm} If $A$ is a commutative $C^*-$algebra, then $A$ and $A^{**}$ are weakly amenable.
\end{proposition}

Let $G$ be a locally compact group with the identity element $e$.  Let us recall that a continuous function $\omega: G\rightarrow [1, \infty)$ is called a \emph{weight} \emph{function} if for every $x, y\in G$
$$
\omega(xy)\leq\omega(x)\;\omega(y)\quad\hbox{and}\quad\omega(e)=1.
$$
Let $C_b (G, 1/\omega)$ denote the space of all  complex-valued functions $f$ on $G$ such that $f/\omega\in C_b(G)$, the space of all bounded continuous functions on $G$. Let also $C_0(G, 1/\omega)$ denote the
subspace of $C_b(G, 1/\omega)$ consisting of all functions  that vanish at
infinity.
We now give another consequence of Proposition \ref{ccmm}.

\begin{corollary} Let $G$ be a locally compact group. The following statements hold.

\emph{(i)} $C_0(G, 1/\omega)$ and $M(G, \omega)^*$ are weakly amenable.

\emph{(ii)} $C_b(G, 1/\omega)$ and $C_b(G, 1/\omega)^{**}$ are weakly amenable.
\end{corollary}

The function $\Omega: G\times G\rightarrow (0, 1]$ defined by
$$
\Omega(x, y)=\frac{\omega(xy)}{\omega(x)\omega(y)}
$$
 is called \emph{zero cluster} if for every pair of sequences $(x_n)_n$ and $(y_m)_m$ of distinct elements in $G$, we have
\begin{eqnarray*}\label{positive cluster}
\lim_n\lim_m\Omega(x_n, y_m)=0=\lim_m\lim_n\Omega(x_n, y_m),
\end{eqnarray*}
whenever both iterated limits exist. A function $f\in C_b (G, 1/\omega)$ is called $\omega-$\emph{weakly almost periodic} (respectively, $\omega-$\emph{almost periodic}) if the set
$$
\{\frac{_xf}{\omega(x)\omega}: x\in G\}
$$
is relatively weakly (respectively, norm) compact in $C_b (G)$, where $_xf(y)=f(yx)$ for all $x, y\in G$. The set of all $\omega-$(respectively, $\omega-$weakly) almost periodic on $G$ is denoted by $\hbox{Ap}(G, 1/\omega)$ (respectively, $\hbox{Wap}(G, 1/\omega)$). The authors \cite{mr} proved that $$\hbox{Wap}(G, 1/\omega) = C_b(G, 1/\omega)$$ if and only if $G$ is compact or $\Omega$ is zero cluster. They also showed that
$\hbox{Ap}(G, 1/\omega)= C_b(G, 1/\omega)$ if and only if $G$ is either compact or discrete and
$\Omega\in C_0(G\times G)$. These facts proved the following result.

\begin{corollary} Let $G$ be a locally compact group. Then the following statements hold.

\emph{(i)} If $G$ is discrete and $\Omega\in C_0(G\times G)$, then $\emph{Ap}(G, 1/\omega)$ is weakly amenable.

\emph{(ii)} If $\Omega$ is zero cluster, then $\emph{Wap}(G, 1/\omega)$ is weakly amenable.

\emph{(iii)} If $G$ is compact, then $\emph{Ap}(G, 1/\omega)$ and $\emph{Wap}(G, 1/\omega)$ are weakly amenable.
\end{corollary}

Assume that $L^1(G, \omega)$ is the Banach space of all Borel measurable functions $f$ on $G$ such that $\omega f\in L^1(G)$, the group  algebra of $G$. Then $L^1(G, \omega)$ with the convolution product ``$\ast$" and the norm $\|f\|_\omega=\|\omega f\|_1$ is a Banach algebras \cite{dl, r0}. Let $L^\infty (G,1/\omega)$ be the Banach algebra of all Borel measurable
functions $f$ on $G$ with $f/\omega\in L^\infty (G)$, the Lebesgue space of bounded Borel measurable functions on $G$.
A Borel measurable complex-valued function $q$ on $G\times G$ is called a \emph{quasi-additive function} if for almost every where $x, y, z\in G$
$$q(xy, z)=q(x, yz)+q(y, zx).$$
If there exists $h\in L^\infty(G, 1/\omega)$ such that $$q(x, y)= h(xy)-h(yx)$$ for almost every where $x, y\in G$, then $q$ is called \emph{inner}. From Theorem 2.4 in \cite{mrl} and Corollary \ref{bil} we have the following.

\begin{theorem}\label{nar1} Let $G$ be a locally compact group. Then the following assertions are equivalent.

\emph{(a)} $L^1(G, \omega)$ is  weakly amenable.

\emph{(b)} $L^1(G, \omega)$ is  cyclically weakly amenable.

\emph{(c)} $L^1(G, \omega)$ is point amenable.

\emph{(d)}  Every non-inner quasi-additive function in $L^\infty(G, 1/\omega)$ is unbounded.
\end{theorem}

Let $CD(G, \omega)$ be the set of all quasi-additive functions $q: G\times G\rightarrow\Bbb{C}$ such that  $q(x, e)=0$ for all almost ever where $x\in G$ and
$$\sup_{x, y\in G}\frac{|q(x, y)|}{\omega(x)\omega(y)}<\infty.$$
Note that for every $x, y\in G$, we have $$q(xy, e)=q(x, y)+q(y, x).$$  This shows that $p\in CD(G, \omega)$ if and only if $$q(x, y)+q(y, x)=0$$ for all almost ever where $x, y\in G$. Also note that  $q(x, x)=(1/2) q(x^2, e)$ for all almost ever where $x\in G$. So if $q\in CD(G, \omega)$, then $q(x, x)=0$.

\begin{theorem}\label{nar} Let $G$ be a locally compact group. Then  $L^1(G, \omega)$ is  cyclically amenable if and only if every element of $CD(G, \omega)$ is inner.
\end{theorem}
{\it Proof.} Let $q\in D(G, \omega)$. By Theorem 2.3  (i) in \cite{mrl}, there exists a continuous derivation
$D: L^1(G, \omega)\rightarrow L^1(G, \omega)^*$ such that
$$
\langle D(f), g\rangle=\langle q, f\otimes g\rangle
$$
for all $f, g\in L^1(G, \omega)$. It follows that
\begin{eqnarray}\label{1256}
\langle D(f), g\rangle+\langle D(g), f\rangle&=&\langle q, f\otimes g\rangle+\langle q, g\otimes f\rangle\nonumber\\
&=&\int_G\int_Gq(x, y)f(x)g(y)\; dxdy\nonumber\\
&+&\int_G\int_Gq(y,x)f(x)g(y)\;dxdy\\
&=&\int_G\int_G(q(x, y)+q(y,x))f(x)g(y)\;dxdy\nonumber
\end{eqnarray}
for all $f,g\in L^1(G, \omega)$.

 Assume now that $q\in CD(G, \omega)$. Then $$q(x, y)+q(y, x)=0$$ for all $x, y\in G$. So by (\ref{1256}), $D$ is cyclic. If $L^1(G, \omega)$ is cyclically amenable, then $D$ is inner. Apply Theorem 2.3 (ii) in \cite{mrl} to conclude that $q$ is inner.

Conversely, let $D: L^1(G, \omega)\rightarrow L^1(G, \omega)^*$ be a cyclic derivation. From (\ref{1256}) we see that $q\in CD(G, \omega)$. By assumption, $q$ is inner and so $D$ is inner. Thus $L^1(G, \omega)$ is  cyclically amenable.$\hfill\square$\\

Let $M(G)$ be the measure algebra of $G$ as defined in \cite{hr}.
In view of Theorem 1.2 in \cite{dgh} and Theorem \ref{sem}, we observe the next result.

\begin{theorem} Let $G$ be a locally compact group. Then the following assertions are equivalent.

\emph{(a)} $M(G)$ is weakly amenable.

\emph{(b)} $M(G)$ is  cyclically weakly amenable.

\emph{(c)} $M(G)$ is point amenable.

\emph{(d)} $G$ is discrete.
\end{theorem}

The next result is proved by using a routine argument. So we omit it.

\begin{proposition} Let $A$ be a Banach algebra generated by the set $B$. Then $A$ is point amenable if and only if every continuous point derivation of $A$ is zero on $B$.
\end{proposition}

\begin{theorem} Let $A$ be a unital Banach algebra, $a\in A$ and $B$ be the Banach algebra generated by the set $\{1, a\}$. Then the following assertions are equivalent.

\emph{(a)} $B$ is weakly amenable.

\emph{(b)} $B$ is cyclically weakly amenable.

\emph{(c)} $B$ is point amenable.

\emph{(d)} For every continuous point derivation $d$ at $\varphi\in\Delta(B)$ of $B$, $d(a)=0$.\\
The Banach algebra $B$ is always cyclically amenable.
\end{theorem}
{\it Proof.} Since $B$ is a unital, commutative Banach algebra, by Theorem \ref{tasisat}, the statements (a), (b) and (c) are equivalent.  Clearly, (c) implies (d). Now, let $\varphi\in\Delta(B)$ and $d$ be a continuous point derivation at $\varphi$ of $B$. Then we have
$$
d(1)=d(1.1)=2d(1)\varphi(1) =2d(1).
$$
So $d(1)=0$.  So (d)$\Rightarrow$(c).

To complete the proof, let $D: B\rightarrow B^*$ be a continuous cyclic derivation. Then for every $m,n\in\Bbb{Z}$ there exists $k\in\Bbb{Z}$ such that
$$
k\langle D(a), a^{n+m-1}\rangle=\langle D(a^n), a^m\rangle+\langle D(a^m), a^n\rangle=0.
$$
So for every $x\in B$, we have $\langle D(a), x\rangle=0$. Since $D$ is cyclic, $\langle D(x), a\rangle=0$. It follows that $\langle D(a^n), a^m\rangle=0$ for all $m,n\in\Bbb{Z}$. Thus $D=0$. From this and commutativity of $B$ we infer that $B$ is cyclically amenable.$\hfill\square$\\

Let $P(x)$ be the Banach algebra of all polynomials in single variable $x$.

\begin{corollary} Then the following assertions are equivalent.

\emph{(a)} $P(x)$ is weakly amenable.

\emph{(b)} $P(x)$ is cyclically weakly amenable.

\emph{(c)} $P(x)$ is point amenable.

\emph{(d)} For every continuous point derivation $d$ of $P(x)$, $d(x)=0$.\\
 The Banach algebra $P(x)$ is always cyclically amenable.
\end{corollary}

One can show that commutative Banach algebras spanned by its idempotents are weakly amenable; see Proposition 2.8.72 in \cite{d}. Here, we investigate this result for point amenability,  cyclic weak amenability and cyclic amenability.

\begin{proposition}\label{idempotent} Let $A$ be a Banach algebra spanned by its idempotents. Then the following statements hold.

\emph{(i)} $A$ is \emph{0}-point amenable.

\emph{(ii)} If $A$ is commutative, then $A$ is cyclically amenable and  cyclically weakly amenable.

\end{proposition}
{\it Proof.} (i) Let $d$ be a continuous  point derivation at $\varphi\in\Delta_0(A)$. If $e\in A$ is idempotent, then
$$
d(e)=d(e^2)=2\varphi(e)d(e).
$$
Since $\varphi(e)=0$ or $\varphi(e)=1$, it follows that $d(e)=0$. Hence
$$
d(ee^\prime)=d(e)\varphi(e^\prime)+\varphi(e)d(e^\prime)=0
$$
for all idempotent elements $e$ and $e^\prime$ in $A$. Continuing by induction and using continuity of $d$ we obtain $d=0$.

(ii) This follows from Proposition 2.8.72 in \cite{d} and Theorem \ref{shah}.$\hfill\square$

\vspace{2mm}

 {\footnotesize
\noindent {\bf Mohammad Javad Mehdipour}\\
Department of Mathematics,\\ Shiraz University of Technology,\\
Shiraz
71555-313, Iran\\ e-mail: mehdipour@.ac.ir\\
{\bf Ali Rejali}\\
Department of Pure Mathematics,\\ Faculty of Mathematics and Statistics,\\ University of Isfahan,\\
Isfahan
81746-73441, Iran\\ e-mail: rejali@sci.ui.ac.ir\\

\footnotesize

\begin{thebibliography}{99}

\bibitem{bcd} W. Bade, P. C. Curtis and H. G. Dales, Amenability and weak amenability for Beurling and Lipschitz algebras, Proc. London Math. Soc., (3) 55 (1987), no. 2, 359--377. 

\bibitem{bd} F. F. Bonsall and J. Duncan, Complete Normed Algebras, Ergebnisse der Mathematik und ihrer Grenzgebiete, Band 80, Springer-Verlag, New York-Heidelberg, 1973.

\bibitem{b} C. R. Borwick, Johnson-Hochschild Cohomology of Weighted Group Algebras and Augmentation Ideals, Ph. D. thesis, University of Newcastle upon Tyne, 2003.

\bibitem{c} J. B. Conway, A Course in Functional Analysis, Second edition, Graduate Texts in Mathematics, 96, Springer-Verlag, New York, 1990.

\bibitem{d} H. G. Dales, Banach Algebras and Automatic Continuity,  London Mathematical Society Monographs, New Series, 24, Oxford Science Publications, The Clarendon Press, Oxford University Press, New York, 2000.

\bibitem{dgh} H. G. Dales, F. Ghahramani and A. Y. A. Helemskii, The amenability of measure algebras, J. London Math. Soc., (2) 66 (2002) 213--226.

\bibitem{dl} H. G. Dales and A. T. Lau, The second duals of Beurling algebras, Mem. Amer. Math. Soc., 177 (2005), no. 836, 1--191. 

\bibitem{f} B. Forrest, Weak amenability and the second dual of the Fourier algebra, Proc. Amer. Math. Soc., 125 (1997), no. 8, 2373--2378.

\bibitem{fm} B. Forrest and L. W. Marcoux, Weak amenability of triangular Banach algebras, Trans. Amer. Math. Soc., 354 (2002), no. 4, 1435--1452.

\bibitem{fr} B. Forrest and V. Runde, Amenability and weak amenability of the Fourier algebra, Math. Z., 250 (2005), no. 4, 731--744.

\bibitem{fss} B. Forrest, E. Samei and N. Spronk, Weak amenability of Fourier algebras on compact groups, Indiana Univ. Math. J., 58 (2009), no. 3, 1379--1393.

\bibitem{gl} F. Ghahramani and J. Laali, Amenability and topological centres of the second duals of Banach algebras, Bull. Austral. Math. Soc., 65 (2002), no. 2, 191--197.

\bibitem{gro3} N. Gronback, Amenability of weighted discrete convolution algebras on cancellative semigroups, Proc. Royal Soc. Edinburgh Sect. A, 110 (1988), no. 3-4, 351--360.

\bibitem{gro} N. Gronback, A characterization of weakly amenable Banach algebras, Studia. Math., 94 (1989), no. 2, 149--162.

\bibitem{gro4} N. Gronback, Commutative Banach algebras, module derivations, and semigroups, J. London Math. Soc., (2) 40 (1989), no. 1, 137--157.

\bibitem{gro1} N. Gronback, Amenability and weak amenability of tensor algebras and algebras of nuclear operators, J. Austral. Math. Soc. Ser. A, 51 (1991), no. 3, 483--488.

\bibitem{gro2} N. Gronback, Weak and cyclic amenability for non-commutative Banach algebras, Proc.
Edinburgh Math. Soc.,  (2) 35 (1992), no. 2, 315--328.

\bibitem{h1} U. Haagerup, All nuclear $C^*-$algebras are amenable, Invent. Math., 74 (1983), no. 2, 305--319.

\bibitem{h2} U. Haagerup and N. J. Laustsen, Weak amenability of $C^*$-algebras and a theorem of Goldstein, Banach algebras '97 (Blaubeuren), 223--243, de Gruyter, Berlin, 1998.

\bibitem{hel1} A. Ya. Helemskii, Some remarks about ideas and results of topological homology,  Conference on Automatic Continuity and Banach Algebras (Canberra, 1989), 203--238, Proc. Centre Math. Anal. Austral. Nat. Univ., 21, Austral. Nat. Univ., Canberra, 1989.

\bibitem{hr} E. Hewitt and K. Ross, Abstract Harmonic Analysis I, Springer-Verlag, New York, 1970.

\bibitem{ho1} G. Hochschild, On the cohomology groups of an associative algebra, Ann. of  Math., (2) 46 (1945), 58--67. 


\bibitem{ho2} G. Hochschild, On the cohomology theory for associative algebras, Ann. of  Math. (2) 47 (1946), 568--579. 


\bibitem{joh} B. E. Johnson, Cohomology in Banach algebra, Mem. Amer. Math. Soc., 127 (1972).

\bibitem{j1} B. E. Johnson, Derivations from $L^1(G)$ into $L^1(G)$ and $L^\infty(G)$, Harmonic analysis (Luxembourg, 1987), 191--198, Lecture Notes in Math., 1359, Springer, Berlin, 1988.

\bibitem{kam} H. Kamowitz, Cohomology groups of commutative Banach algebras, Trans. Amer.
Math. Soc., 102 (1962), 352--372. 


\bibitem{kams} H. Kamowitz and S. Scheinberg, Derivations and automorphisms of $L^1(0, 1)$, Trans. Amer. Math. Soc., 135 (1969), 415--427.

\bibitem{mr} M. J. Mehdipour and A. Rejali, Regularity and amenability of weighted Banach algebras and their second dual on locally compact groups,  arXiv:2112.13286v1.

\bibitem{mrl} M. J. Mehdipour and A. Rejali, Weak amenability of weighted group algebras, arXiv:2209.08346.

\bibitem{p2} A. Pourabbas, Second cohomology of Beurling algebras, Saitama Math. J., 17 (1999), 87--94.

\bibitem{r0} H. Reiter and J. D. Stegeman, Classical Harmonic Analysis and Locally Compact Groups, Second edition, London Mathematical Society Monographs, New Series, 22, The Clarendon Press, Oxford University Press, New York, 2000.

\bibitem{s2} D. R. Sherbert, The structure of ideals and point derivations in Banach algebras of Lipschitz functions, Trans. Amer. Math. Soc., 111 (1964), 240--272.

\bibitem{s3}  D. R. Sherbert, Banach algebras of Lipschitz functions, Pacific J. Math., 13 (1963), 1387--1399.

\bibitem{s4}  D. R. Sherbert, Banach Algebras of Lipschitz Functions, Thesis (Ph.D.)–Stanford University, ProQuest LLC, Ann Arbor, MI, 1962. 102.
\end{thebibliography}
\end{document}